\documentclass[12pt]{amsart}
\usepackage{latexsym, amsmath, amscd, amssymb, amsthm} 
\newtheorem*{te}{Theorem}
\newtheorem{lm}{Lemma}

\begin{document}

 \title{   On complete system of invariants for the   binary form of  degree 7} 

\author{Leonid Bedratyuk} \address{ Khmel'nyts'ky National University, Instytuts'ka st. 11, Khmel'nyts'ky , 29016, Ukraine}
\email {bedratyuk@ief.tup.km.ua}
\begin{abstract}
A minimal system of homogeneous generating elements of  the invariants algebra for  the  binary form of degree 7 is calculated.
\end{abstract}

\maketitle

\section{Introduction.}

Let $V_n$ be a vector $\mathbb{C}$-space of the binary forms of degree $d$ considered with natural action of the group  $G=SL(2,\mathbb{C}).$ Let us extend the action of the group  $G$ to the polynomial functions algebras $\mathbb{C}[V_d]$ and $\mathbb{C}[V_d \oplus \mathbb{C}^2 ].$ 
Denote by  $I_d=\mathbb{C}[V_d]^{\,G}$ and by  ${C_d=\mathbb{C}[V_d \oplus \mathbb{C}^2 ]^{\,G}}$ the corresponding subalgebras of   $G$-invariant functions. 
In the vocabulary of classical invariant theory  the algebras  $I_d$ and   $C_d$ are called,  respectively,    algebra invariants and  algebra covariants of the  binary form of d-th degree.
Let  $I_d^{+}$ be an ideal of  $I_d$ generated by all homogeneus elements of positive power. Denote by  $\bar I_d$ a set of  homogeneus elements of $I_d^{+}$ such that their images in  $I_d^{+}/(I_d^{+})^2$ form a basis of the vector space. The  set  $\bar I_d$ is caled complete system of invariants of the  d-th degree binary form. Elements of   $\bar I_d$ form a minimal system of homogeneous generating elements of the invariants algebra $I_d.$ Denote by $n_d$ a number of  elements of the set    $\bar I_d.$ In the same way we may determine a   complete system $\overline C_d$ of coinvariants of the   binary form of the degree  $d$

The complete systems of invariants and covariants was a topic of major research interest in classical invariant theory of the 19th century.
It is easy to show that $n_1=0, $  $n_2=1,$ $n_3=1.$ A complete system of invariants and covariants in the case   $d=4$ was calculated by Bool, Cayley, Eisenstein and in the case $d=5$ by common efforts of Cayley and Hermite, see survey \cite{Dix}.
The complete systems  of invariants and covariants in the case   $d=6$ were calculated by Gordan, see  \cite{Gor}. In particular, $n_4=2,$  $n_5=4,$ $n_6=5.$ A complete system of 9 invariants in the case  $d=8$ were computed by  Gall  \cite{Gall} and Shioda \cite{Shi}.

The case $d=7$ was considered by  Sylvester and Gall. By using Cayley-Sylvester's technique Sylvester \cite{Sylv}-\cite{SF} found a low bound  for the invariants number, namely $n_7\geq  28.$ In the  paper   \cite{Gall-1} Gall essentially   showed  that  $n_7\leq 33.$ Nearly 100 years later  Dixmier and  Lazard  \cite{DL} proved  that  $n_7= 30.$ However, the question -- what elements exactly form the complete system of invariants remained open until recently. Gall found the 33 invariants in an implicity way using the sybolic method.
The invariant's symbolic representation though simple is, nevertheless, much too hard to check for irreducibility of the invariant. The verification could be done for  invariants in their explicit representation. Unfortunately, some of those offerred by Gall infeasible for computer calculation due to their high tranvectant's order.

To solve the computation problem we offer a form of represenation of the invariants, which  is an intermediate form between the highly unwieldy explicit represenation and too "compressed" symbolic represenation. A first step in the simplifycation of a calculation is  calculation  semi-invariants instead calculation of covariants. Let us consider a covariant  as polynomial  of generating  functions of the polynomial function algebra   $\mathbb{C}[V_d \oplus \mathbb{C}^2 ].$ Then a semi-invariat is just a leading coefficient of the polynomial with respect to usual lexicographical ordering. A semi-invariant is an invariant of upper unipotent matrix  subalgebra of Lie algebra $\mathfrak{sl_{2}}.$ 

Let us indentify the algebra $\mathbb{C}[V_d]$  with the algebra $\mathbb{C}[X_d]:=\mathbb{C}[t,x_1,x_2,\ldots, x_d],$ and the algebra $\mathbb{C}[V_d \oplus \mathbb{C}^2 ]$   with the polynomial algebra $\mathbb{C}[t,x_1,x_2,\ldots, x_n,Y_1,Y_2].$

The generating elements  $\Bigl( \begin{array}{ll}  0\, 1 \\ 0\,0 \end{array} \Bigr),$ $\Bigl( \begin{array}{ll}  0\, 0 \\ 1\,0 \end{array} \Bigr)$ of the tangent Lie algebra  $\mathfrak{sl_{2}}$ act    on  $\mathbb{C}[V_d]$  by  derivations
$$
\begin{array}{l}
D_1:=\displaystyle t\frac{\partial}{\partial x_1}+2\, x_1\frac{\partial}{\partial x_2}+\cdots +n\,x_{d-1}\frac{\partial}{\partial x_d}, \\
D_2:=\displaystyle d\,x_1\frac{\partial}{\partial t}+(d-1)\,x_2\frac{\partial}{\partial x_1}+\cdots +x_{d}\frac{\partial}{\partial x_{d-1}}.
\end{array}
$$

It follows that the invariant algebra $I_d$  coincides with an algebra of polynomial solutions of the following first order PDE system, see \cite{Hilb}, \cite{Gle}:

$$
\left\{
\begin{array}{r}
\displaystyle t\frac{\partial u}{\partial x_1}+2\, x_1\frac{\partial u}{\partial x_2}+\cdots +n\,x_{d-1}\frac{\partial u}{\partial x_d}=0, \\
\displaystyle d\,x_1\frac{\partial u}{\partial t}+(d-1)\,x_2\frac{\partial u}{\partial x_1}+\cdots +x_{d}\frac{\partial u}{\partial x_{d-1}}=0,
\end{array}
\right.
\eqno (*)
$$
i.e.  $I_d=\mathbb{C}[X_d]^{D_1}  \cap \mathbb{C}[X_d] ^{D_2} ,$ where $u \in \mathbb{C}[X_d],$  and $$\mathbb{C}[X_d]^{D_i}:=\{ f \in \mathbb{C}[X_d]| D_i(f)=0 \}, i=1,2.$$

The algebra of semi-invariants coincides with the algebra  $\mathbb{C}[X_d]^{D_1}.$ It is easy to get an explicit form of the  algebra $\mathbb{C}[X_d]^{D_1},$ see, for example, \cite{Ess}. Namely -- $$\mathbb{C}[X_d]^{D_1}=\mathbb{C}[t, z_2,\ldots, z_d][\frac{1}{t}] \cap \mathbb{C}[X_d],$$ here  $z_i$ are some  functional independent semi-invariants of degree  $i.$ The polynomials $z_i$ arised in the first time in Cayley, see  \cite{Gle}. Therefore,  any semi-invariant we may write as rational fraction of ${\displaystyle \mathbb{C}[Z_d][\frac{1}{t}]:=\mathbb{C}[t, z_2,\ldots, z_d][\frac{1}{t}].}$  This form of semi-invariants is more compact than their standard form as a polynomial of   $\mathbb{C}[X_d].$ In this case,
 a semi-invariant's number of terms is less by  tens times  than  the number of terms of the corresponding covariant. Thus making a substantial difference computationally.  
 By using Robert's theorem \cite{Rob}, knowing a semi-invariant one may restore the corresponding  covariant.  

For computing of  the semi-invariants we introduce  a semitransvectant that is an analogue of a transvectant. In the paper we have found effective formula for the semitransvectants computation and for the case $d=7$ and have calculated all irreducible semi-invariants up to 13 degree.  An invariant of degree $n$ we are looking as a semitranvectant  of the form  $[u,v]^i$, where  $u, v$ are semi-invariants of the degree $n/2.$ Therefore, by using the obtained semi-invariants  we have computing all invariants up to 26 degree. An  invariant of the degree 30 was  taken  from Gall's paper  \cite{Gall}. We compute the invariant and check out that the invariant is irreducible. Thus, in the paper in explicit way were computed complete system of 30 invariants of the  binary form of degree 7. 

All calculation were done with Maple.

%-----------------------------------------------------------------------------------------

\section{Premilinaries.}

%--------------------------------------------------------------------------------------------------------------

Before any calculation, we first try to simplify both the representations and calculations of covariants. 
Let  ${\it \kappa} : C_d \longrightarrow \mathbb{C}[X_d]^{D_1}$ be the   $\mathbb{C}$-linear map takes each homogeneous covariant of order $k$ to his leading coefficient, i.e. a coefficient of   $Y_1^k.$ Follow by classical  tradition an element of the algebra   $\mathbb{C}[X_d]^{D_1}$  is called  {\it semi-invariants, } a degree of a homogeneous covariant with respect to the variables set $X_d$ is called {\it degree} of the covariant and its  degree with respect to the variables set $Y_1, Y_2$  is called {\it order.}

Suppose  $F=\displaystyle \sum_{i=0}^m \, f_i { m \choose i } Y_1^{m{-}i} Y_2^i$ be a covariant of order $m,$  $\kappa(F)=f_0 \in \mathbb{C}[X_d]^{D_1}.$ The classical Robert's theorem \cite{Rob}  states that the covariant $F$ is completelly and   uniquely determined  by its leading coefficient $f_0,$ namely 
$$
F=\sum_{i=0}^{m} \frac{D_2^i(f_0)}{i!} Y_1^{m-i}Y_2^i.
$$
On the other hand, every semi-invariant is a leading coefficient of some covariant, see \cite{Gle}. This 
give us well defined   explicit form of the inverse map  ${\kappa^{-1} :   \mathbb{C}[X]^{d_1} \longrightarrow C_d,}$ namely

$$
\kappa^{-1}(a)=\sum_{i=0}^{{\rm ord}(a)} \frac{D_2^i(a)}{i!} Y_1^{{\rm ord}(a)-i}Y_2^i,  
$$
here  $ a \in  \mathbb{C}[X]^{d_1}$ and  ${\rm ord}(a)$ is an order of the element $a$ with respect to the locally nilpotent derivation  $D_2,$  i.e. ${\rm ord}(a)\!:=\max \{ s, D_2^s(a) \neq 0 \}.$  It is evident that for any $a \in  \mathbb{C}[X]^{d_1}$ the order ${\rm ord}(a)$  coincides with the  order of the corresponding covariant $\kappa^{-1}(a).$ For example, since ${\rm ord}(t)=d,$ we have 
$$
\kappa^{-1}(t)=\sum_{i=0}^{{\rm ord}(t)} \frac{D_2^i(t)}{i!} Y_1^{{\rm ord}(t)-i}Y_2^i =t Y_1^d+\sum_{i=1}^{d} { d \choose i } x_i Y_1^{d-i}Y_2^i.
$$
As we see  the $\kappa^{-1}(t)$ is just the  basic binary form of order $d.$ From polynomial functions point of view the covariant $\kappa^{-1}(t)$ is the evaluation  map.

It is clear that an invariant is a semi-invariant of the zero order. 
Thus,  the  problem of finding of a complete system of the  algebra   $\overline C_d$ is equivalent to the  problem  of finding of a complete system of semi-covariants's algebra   $\mathbb{C}[X]^{D_1}.$ For one's turn, the  problem of finding of complete system of the  algebra $\bar I_d$  is equivalent to the  problem  of finding of a subsystem in  $\overline C_d$ that is generated by  elements of zero order. It are well known classical results.

A structure of  constants algebras for locally nilpotent  derivations can be easy determined, see for example  \cite{Ess}.  In particular, for the derivation $D_1$  we get 
$$
\mathbb{C}[X_d]^{D_1}=\mathbb{C}[t,\sigma(x_2),\dots ,\sigma(x_n)][\frac{1}{t}] \cap \mathbb{C}[X_d],
$$
\noindent
where   $\sigma: \mathbb{C}[X_d] \to \mathbb{C}(X_d)^{D_1} $  is a ring homomorphism defined by  
$$
\sigma(a)=\sum_{i=0}^{\infty} d_1^{\,i}(a) \frac{\lambda^i}{i!}, \lambda = -\frac{x_1}{t}.
$$
After not complicated simplification we obtain  $\displaystyle \sigma(x_i)=\frac{z_{i+1}}{t^i},$ where  $z_i \in  \mathbb{C}[X_d]^{D_1}$ and 
$$
z_i:= \sum_{k=0}^{i-2} (-1)^k {i \choose k} x_{i-k}  x_1^k t^{i-k-1} +(i-1)(-1)^{i+1} x_1^i, i=2,\ldots,n.
$$
Especially
$$
\begin{array}{l}
z_2={x_{2}}\,t - {x_{1}}^{2}
\\
z_3={x_{3}}\,t^{2} + 2\,{x_{1}}^{3} - 3\,{x_{1}}\,{x_{2}}\,t
\\
z_4={x_{4}}\,t^{3} - 3\,{x_{1}}^{4} + 6\,{x_{1}}^{2}\,{x_{2}}\,t - 4
\,{x_{1}}\,{x_{3}}\,t^{2}
\\
z_5={x_{5}}\,t^{4} + 4\,{x_{1}}^{5} - 10\,{x_{1}}^{3}\,{x_{2}}\,t + 
10\,{x_{1}}^{2}\,{x_{3}}\,t^{2} - 5\,{x_{1}}\,{x_{4}}\,t^{3}
\\
z_6={x_{6}}\,t^{5} - 5\,{x_{1}}^{6} + 15\,{x_{1}}^{4}\,{x_{2}}\,t - 
20\,{x_{1}}^{3}\,{x_{3}}\,t^{2} + 15\,{x_{1}}^{2}\,{x_{4}}\,t^{3}
 - 6\,{x_{1}}\,{x_{5}}\,t^{4}
\\
z_7={x_{7}}\,t^{6} + 6\,{x_{1}}^{7} - 21\,{x_{1}}^{5}\,{x_{2}}\,t + 
35\,{x_{1}}^{4}\,{x_{3}}\,t^{2} - 35\,{x_{1}}^{3}\,{x_{4}}\,t^{3}
 + 21\,{x_{1}}^{2}\,{x_{5}}\,t^{4} - 7\,{x_{1}}\,{x_{6}}\,t^{5}
\end{array}
$$
Thus we obtain   
$$
\mathbb{C}[X_d]^{D_1}=\mathbb{C}[t,z_2,\ldots,z_d][\frac{1}{t}] \cap \mathbb{C}[X_d].
$$
Hence, a generating elements of the semi-invariant algebra  $\mathbb{C}[X_d]^{D_1}$ we may consider  as a rational fraction $\displaystyle \frac{f(z_2,\ldots,z_n)}{t^s},$ $f \in \mathbb{C}[Z_d]:=\mathbb{C}[t,z_2,\ldots,z_d],$ $s \in  \mathbb{Z_{+}}.$

To make  a calculation with an invariants in such representation we need to know an  action of the operator  $D_2$ in new coordinates $t,z_2,\ldots,z_d.$
  Denote by  $D$  extention  of the derivation  $D_2$  to the algebra $\displaystyle \mathbb{C}[Z_d][\frac{1}{t}]:$
$$
D:=D_2(t)\,\frac{\partial}{\partial t}+D_2(z_2)\,\frac{\partial}{\partial z_2}+\ldots +D_2(z_n)\,\frac{\partial}{\partial z_n}.
$$ 
In  \cite{Aut} the autor proved  that 
$$
\begin{array}{l}
D(t)=-n t \lambda,\\

D(\sigma(x_2))=\displaystyle (n-2)\sigma(x_3)-(n-4)\sigma(x_2) \lambda,\\
D(\sigma(x_i))=\displaystyle (n-i)\sigma(x_{i+1})-(n-2i)\sigma(x_i) \lambda-i(n-1) \frac{\sigma(x_2) \sigma(x_{i-1})}{t}, \mbox{for  }  i>2.
\end{array}
$$
 Taking into account   $\displaystyle \sigma(x_i)=\frac{z_{i+1}}{t^i},$ $ \displaystyle \lambda = -\frac{x_1}{t}$  we can obtain an expression and for  $D(z_i),$ ${ i=2,\ldots,d.}$ Especially, for  $d=7$ we get :
$$ 
\begin{array}{l}

D=\displaystyle 7\,{x_{1}}\,{\frac {\partial }{\partial t}} - 
\displaystyle \frac {( - 15\,{x_{1}}\,{z_{3}} + 18\,{z_{2}}^{2}
 - 4\,{z_{4}})}{t}{\frac {\partial }{\partial {z_{3}}}} 
 + {\displaystyle \frac {(20\,{x_{1}}\,{z_{4}} - 24\,{z_{2}}\,{z
_{3}} + 3\,{z_{5}})}{
t}}\,{\frac {\partial }{\partial {z_{4}}}} +  \\
\\
\mbox{} + {\displaystyle \frac {(2\,{z_{6}} + 25\,{x_{1}}\,{z_{5}
} - \displaystyle 30\,{z_{2}}\,{z_{4}})}{t}}\,{\displaystyle \frac {\partial }{\partial {z_{5}}}}
  + {\displaystyle \frac {({z_{7}} + 30\,{x_{1}}\,{z_{6}
} - 36\,{z_{2}}\,{z_{5}})}{t}} \,{\displaystyle \frac {\partial }{\partial {z_{6}}}}+
 \\
\\
\mbox{} + {\displaystyle \frac {7\,(5\,{x_{1}}\,{z_{7}} - 6\,{z_{
2}}\,{z_{6}})}{t}} \,{\displaystyle \frac {\partial }{\partial {z_{7}}}}
 + {\displaystyle \frac {5\,(2\,{x_{1}}\,{z_{2}} + {z_{3}})}{t}} \,
{\displaystyle \frac {\partial }{ \partial {z_{2}}}}.
\end{array}
$$
To calculate the semi-invariants we need have an analogue of the transvectants. 
Suppose $$F=\sum_{i=0}^m \, f_i { m \choose i } Y_1^{m{-}i} Y_2^i, \mbox{   }  G=\sum_{i=0}^k \, f_i { k \choose i } Y_1^{k{-}i} Y_2^i,   \mbox{   }        f _i, g_i \in \mathbb{C}[Z_d][\frac{1}{t}], $$ are two covariants of the orders  $m$ and $ k$ respectively. Let 
$$
(F,G)^r=\sum_{i=0}^r (-1)^i { r \choose i } \frac{\partial^r F}{\partial Y_1^{r-i} \partial Y_2^i}   \frac{\partial^r G}{\partial Y_1^{i} \partial Y_2^{r-i}},
$$
be their   $r$-th transvectant. 
The following lemma give us rule how to find the semi-invariant $\kappa( (F,G)^r)$ without of direct  computing of the covariant  $(F,G)^r.$
\begin{lm}
The leading coefficient  $\kappa ((F,G)^r)$ of the covariant $(F,G)^r,  0 \leq r \leq \min(m,k)$ is calculating by the formula
$$
\kappa((F,G)^r)=\sum_{i=0}^r (-1)^i { r \choose i } \frac{D^i(\kappa(F))}{[m]_i} {\Bigl | _{x_1=0,\ldots,x_r=0}}  \frac{D^{r-i}(\kappa(G))}{[k]_{r-i}}{\Bigl | _{x_1=0,\ldots,x_r=0}},
$$ 
here  $[a]_i:=a (a-1) \ldots (a-(i-1)).$  
\end{lm}
\begin{proof}
In  \cite{Hilb}, p. 87, one may find that  
$$
\kappa((F,G)^r)=\sum_{i=0}^r (-1)^i { r \choose i }f_i g_{r-i}.
$$
By comparing the two different form of the covariant $F$ 
$$
F=\sum_{i=0}^m \, f_i { m \choose i } Y_1^{m{-}i} Y_2^i, \mbox{   and     } F = \sum_{i=0}^{m} \frac{D_2^i(f_0)}{i!} Y_1^{m-i}Y_2^i,
$$
we get  $\displaystyle f_i { m \choose i }= \frac{D_2^i(f_0)}{i!},$  and  $\displaystyle f_i=\frac{D^i(f_0)}{[m]_i}=\frac{D^i(\kappa(F))}{[m]_i}.$ Similarly $\displaystyle g_i=\frac{D^i(\kappa(G))}{[k]_i}.$ Therefore,
$$
\kappa((F,G)^r)=\sum_{i=0}^r (-1)^i { r \choose i } \frac{D^i(\kappa(F))}{[m]_i}   \frac{D^{r-i}(\kappa(G))}{[k]_{r-i}}.
$$
The derivation  $D$ acts in such a way that   ${\displaystyle D(\mathbb{C}[Z_d][\frac{1}{t}]) \subset \mathbb{C}[Z_d,x_1][\frac{1}{t}]}.$ Further we have $${\displaystyle D^2(\mathbb{C}[Z_d][\frac{1}{t}]) \subset \mathbb{C}[Z_d,x_1,x_2][\frac{1}{t}]},$$ and for arbitrary $r$ we get  $${\displaystyle D^r(\mathbb{C}[Z_d][\frac{1}{t}]) \subset \mathbb{C}[Z_d,x_1,x_2,\ldots, x_r][\frac{1}{t}]}.$$  Therefore and  $$\kappa((F,G)^r) \subset \mathbb{C}[Z_d,x_1,x_2,\ldots, x_r][\frac{1}{t}].$$ On the other hand, $\kappa((F,G)^r)$ allways is a semi-invariant and then must be true the  inclusion   ${\displaystyle \kappa((F,G)^r) \subset \mathbb{C}[Z_d][\frac{1}{t}].}$ Thus, in the expression for  $\kappa((F,G)^r)$ after cancelation all coefficients of   $x_1,$ $x_2,\ldots,$ $x_r$ must be equal to zero.  Hence, 
$$
\kappa((F,G)^r) =\kappa((F,G)^r) {\Bigl | _{x_1=0,\ldots,x_r=0}}=\sum_{i=0}^r (-1)^i { r \choose i } \frac{D^i(\kappa(F))}{[m]_i}   \frac{D^{r-i}(\kappa(G))}{[k]_{r-i}}{\Bigl | _{x_1=0,\ldots,x_r=0}}=
$$
$$
=\sum_{i=0}^r (-1)^i { r \choose i } \frac{D^i(\kappa(F))}{[m]_i} {\Bigl | _{x_1=0,\ldots,x_r=0}}  \frac{D^{r-i}(\kappa(G))}{[k]_{r-i}}{\Bigl | _{x_1=0,\ldots,x_r=0}}.
$$
\end{proof} 
Let   $f,g$ be two semi-invariants. Their numers are a polynomials of  $z_2,\ldots,z_n$ with rational coefficients. Then the semi-invariant  $\kappa((\kappa^{-1}(f),\kappa^{-1}(g))^i)$  be a fraction and their numer be a polynomial  of  $z_2,\ldots,z_n$ with rational coefficients too.  Therefore we may multiply $\kappa((\kappa^{-1}(f),\kappa^{-1}(g))^r)$ by some rational number  $q_r(f,g) \in \mathbb{Q}$ such that the numer of the expression  $ q_r(f,g) \kappa((\kappa^{-1}(f),\kappa^{-1}(g))^r)$  be now a polynomial with an integer coprime coefficients. Put 
$$
[f,g]^r:=q_r(f,g) \kappa((\kappa^{-1}(f),\kappa^{-1}(g))^r), 0 \leq r\leq \min({\rm ord}(f),{\rm ord}(g)). 
$$
The expression $[f,g]^r$ is said to be the  r-th {\it semitranvectant} of the semi-invariants $f$ and $g.$

The following statements are direct consequences of the corresponding  transvectant properties, see \cite{Gle}:
\begin{lm} Let  $f, g$ be two semi-invariants.  Then the folloving conditions hold:
\begin{enumerate}
\item[({\it i})] the semitransvectant $[t, f\,g]^i$ is reducible for  $ 0 \leq i \leq \min(d, \max( {\rm ord}(f),{\rm ord}(g));$
\item[({\it ii})] if  $ {\rm ord}(f)=0,$ then   $[t, f\,g]^i=f [t,g]^i;$
\item[({\it iii})] ${\rm ord}([f,g]^i)={\rm ord}(f)+{\rm ord}(g)-2\, i;$ 
\item[({\it iv})] ${\rm ord}(z_2^{i_1}z_3^{i_3} \cdots z_d^{i_d})=d\,(i_2+i_3+\cdots +i_d)-2\,(2\, i_2+3\, i_3+\cdots +d \, i_d).$
\end{enumerate}
\end{lm}

Let us consider a sample. Suppose $d=7,$ $f=z_2, g=z_3.$ Then  ${\rm ord}(z_2)=7\cdot 2-2\cdot 2=10,$  ${\rm ord}(z_3)=15.$ We have
$$
\begin{array}{l}
d(z_2)={\displaystyle \frac {5\,(2\,{x_{1}}\,{z_{2}} + {z_{3}})}{t}}, [10]_0=1, [10]_1=1,\\ d^{\,2}(z_2)={\displaystyle \frac {10\,(3\,{x_{1}}^{2}\,{z_{2}} + 9\,{x_{1}}\,
{z_{3}} + 6\,{x_{2}}\,{z_{2}}\,t - 9\,{z_{2}}^{2} + 2\,{z_{4}})}{
t^{2}}},[10]_2=10\\
d(z_3)={\displaystyle \frac {15\,{x_{1}}\,{z_{3}} - 18\,{z_{2}}^{2} + 4
\,{z_{4}}}{t}}, [15]_0=1, [15]_1=1,  \\ d^{\,2}(z_3)={\displaystyle \frac {2\,(60\,{x_{1}}^{2}\,{z_{3}} - 252\,{x_{1}}
\,{z_{2}}^{2} + 56\,{x_{1}}\,{z_{4}} + 45\,{x_{2}}\,{z_{3}}\,t - 
138\,{z_{2}}\,{z_{3}} + 6\,{z_{5}})}{t^{2}}}, [15]_2=15,
\end{array}
$$
and
$$
\begin{array}{ll}
d(z_2){\Bigl | _{x_1=0,x_2=0}} ={\displaystyle \frac {5\,( {z_{3}})}{t}},& d^{\,2}(z_2){\Bigl | _{x_1=0,x_2=0}}={\displaystyle \frac {10\,( - 9\,{z_{2}}^{2} + 2\,{z_{4}})}{
t^{2}}},\\
d(z_3){\Bigl | _{x_1=0,x_2=0}}={\displaystyle \frac { - 18\,{z_{2}}^{2} + 4
\,{z_{4}}}{t}}, & d^{\,2}(z_3){\Bigl | _{x_1=0,x_2=0}}={\displaystyle \frac {2\,( - 
138\,{z_{2}}\,{z_{3}} + 6\,{z_{5}})}{t^{2}}}.
\end{array}
$$
We get
$$
\kappa((\kappa^{-1}(z_2),\kappa^{-1}(z_3))^2)=\sum_{r=0}^2 (-1)^r { 2\choose r } \frac{d^r(z_2)}{[10]_r} {\Bigl | _{x_1=0,x_2=0}}  \frac{d^{2-r}(z_3)}{[15]_{2-r}}  {\Bigl | _{x_1=0,x_2=0}}.
$$
After simplifycation we obtain 
$$
\kappa((\kappa^{-1}(z_2),\kappa^{-1}(z_3))^2)=- {\displaystyle \frac {2}{315}} \,{\displaystyle \frac {3\,{z_{
2}}^{2}\,{z_{3}} - 9\,{z_{2}}\,{z_{5}} + 7\,{z_{3}}\,{z_{4}}}{t^{
2}}}.
$$
It is obviousl that  $q_2(z_2,z_3)=- {\displaystyle \frac {315}{2}}.$  Therefore,
$$
\begin{array}{l}
[z_2,z_3]^2={\displaystyle \frac {3\,{z_{
2}}^{2}\,{z_{3}} - 9\,{z_{2}}\,{z_{5}} + 7\,{z_{3}}\,{z_{4}}}{t^{
2}}}= 
- 31\,{x_{1}}^{3}\,{x_{4}}\,t + 16\,{x_{1}}^{4}\,{x_{3}} + 9\,{x
_{1}}^{2}\,{x_{5}}\,t^{2} +\\
+ 7\,{x_{3
}}\,t^{3}\,{x_{4}} +

 30\,{x_{3}}\,t\,{x_{1}}^{2}\,{x_{2}} + 24\,{
x_{1}}\,{x_{2}}\,t^{2}\,{x_{4}} 
\mbox{} - 28\,{x_{3}}^{2}\,t^{2}\,{x_{1}} - 12\,{x_{1}}^{3}\,{x_{
2}}^{2} + 3\,{x_{2}}^{2}\,t^{2}\,{x_{3}}- \\
 - 9\,{x_{2}}\,t^{3}\,{x_{5}} - 9\,{x_{2}}^{3}\,t\,{x
_{1}}.
\end{array}
$$
In the same way 
$$
[t, [z_2,z_3]^2]={\displaystyle \frac {27\,{z_{2}}^{4} - 78\,{z_{2}}^{2}\,{z_{4}}
 - 14\,{z_{4}}^{2} + 69\,{z_{2}}\,{z_{3}}^{2} + 12\,{z_{3}}\,{z_{
5}} + 9\,{z_{2}}\,{z_{6}}}{t^{2}}}=
$$
$$
\begin{array}{l}
= - 78\,{x_{2}}^{2}\,t^{3}\,{x_{4}} - 126\,{x_{2}}^{2}\,{x_{1}}^{4
} + 45\,{x_{2}}^{3}\,t\,{x_{1}}^{2} + 12\,{x_{3}}\,t^{4}\,{x_{5}}
 + 69\,{x_{2}}\,t^{3}\,{x_{3}}^{2} - 173\,{x_{1}}^{2}\,{x_{3}}^{2
}\,t^{2} +\\
\mbox{} + 168\,{x_{1}}^{5}\,{x_{3}} - 249\,{x_{1}}^{4}\,{x_{4}}\,
t{-}9\,{x_{1}}^{2}\,{x_{6}}\,t^{3} + 9\,{x_{2}}\,t^{4}\,{x_{6}}
 + 78\,{x_{1}}^{3}\,{x_{5}}\,t^{2} + 27\,{x_{2}}^{4}\,t^{2}{-}14\,{x_{4}}^{2}\,t^{4}-\\
\mbox{} - 102\,{x_{2}}^{2}\,t^{2}\,{x_{1}}\,{x_{3}} + 303\,{x_{1}
}^{2}\,{x_{2}}\,t^{2}\,{x_{4}} + 78\,{x_{1}}^{3}\,{x_{2}}\,t\,{x
_{3}} + 52\,{x_{4}}\,t^{3}\,{x_{1}}\,{x_{3}} - 90\,{x_{1}}\,{x_{2
}}\,t^{3}\,{x_{5}}.
\end{array}
$$
As we may see from the samples, a representation of semi-invariants as elements of the algebra  
  $\displaystyle \mathbb{C}[Z_d][\frac{1}{t}]$  is more compact than their standard representation as elements of the algebra  $\mathbb{C}[X_d].$ Very rough empire estimate is  that a semi-invariant $\kappa(F)$ has terms number  in ${\displaystyle ([\deg (F)/d]+2)\, ({\rm ord}(F)}+1)$ times less  than the corresponding covariant $F.$ Moreover, from computing  point of view, the semitransvectant formula is more effective than the transvectant  formula. These two favorable circumstances coupled with great Maple powers allow us  to compute  a complete system of invariants of the  7-th degree binary form.

\section{ Computation of an auxiliary semi-invariants. }

Before any invariant's calculation we find all irreducible covariants up to 13th degree inclusive. To do it we use an analogue of the well known $\Omega$-process.
Let  $\overline C_{7,\,i}$ be a subset of elements of  $\overline C_{7}$ that has degree  $i.$ Elements of the set  $\overline C_{7,\,i+1} $ we are seeking  as an irreducible elements of a  basis of a vector space  generated by semitransvectants of the  form  $[t,u\,v]^r, u \in C_{7,\,l},$  $v \in C_{7,\,k},$ $l+k=i,$ $ \max({\rm ord}(u),{\rm ord}(v)) \leq r \leq 7.$ It is a standard linear algebra problem.

The unique semi-invariant of the degree one obviously is  $t, {\rm ord}(t)=7.$
The semitransvectants  $[t,t]^i$ are equal to zero for odd  $i.$ Put 
$$
\begin{array}{ll}
dv_1:=[t,t]^4={\displaystyle \frac {3\,{z_{2}}^{2} + {z_{4}}}{t^{2}}}={x_{4}}\,t - 4\,{x_{1}}\,{x_{3}} + 3\,{x_{2}}^{2} , & {\rm ord}(dv_1){=}6,  \\ dv_2:=[t,t]^6={\displaystyle \frac {{z_{6}} + 15\,{z_{2}}\,{z_{4}} - 10\,{z_{3}, }^{2}}{t^{4}}}={x_{6}}\,t - 6\,{x_{1}}\,{x_{5}} + 15\,{x_{2}}\,{x_{4}} - 10\,{x
_{3}}^{2}, & {\rm ord}(dv_2){=}2, \\
dv_3:=[t,t]^2=z_2={x_{2}}\,t - {x_{1}}^{2}, &  {\rm ord}(dv_2){=}10.
\end{array}
$$
The polynomials     $t^2,$ $dv_1,$ $dv_2,$ $dv_3$   are linear independent. Therefore,   the set  $\overline C_{7,\,2}, $ consists of the irreducible semi-invariants  $dv_1,$ $dv_2,$ $dv_3$ of 2th  degree.

To define the set $\overline C_{7,\,3}$ let us consider the following  19 polynomial $t^3,$ $t \,dv_1,$ $t \,dv_2,$ $t \,dv_3,$ $[t,dv_1]^i,$ $i=1\ldots 6,$  $[t,dv_2]^i,$ $i=1,2$  $[t,dv_3]^i,$ $i=1\ldots 7.$  By direct calculation we select the 6 linearly independent irreducible covariants of 3th degree :
$$
\begin{array}{llll}
tr_1=[t,dv_1]^4, & {\rm ord}(tr_1)=5, & tr_2=[t,dv_3],  & {\rm ord}(tr_2)=15,\\

tr_3=[t,dv_3]^3, &  {\rm ord}(tr_3)=11, & tr_4=[t,dv_3]^4, & {\rm ord}(tr_4)=9,\\

tr_5=[t,dv_3]^5,  & {\rm ord}(tr_5)=7, & tr_6=[t,dv_3]^7, & {\rm ord}(tr_6)=3.

\end{array}
$$

In the same way we may calculate all sets $\overline C_{7,\,i},$ $i\leq 13.$ Let us present the lists of the generating elements.

The set $\overline C_{7,\,4}$ consists of the following 8 ireducible semi-invariants: 
$$
\begin{array}{llll}
ch_1=[t,tr_5]^7, & {\rm ord}(ch_1)=0, & ch_2=[t,tr_3]^7,  & {\rm ord}(ch_2)=4,\\

ch_3=[t,tr_3]^2, &  {\rm ord}(ch_3)=14, & ch_4=[t,tr_3]^4,  & {\rm ord}(ch_4)=10,\\

ch_5=[t,tr_3]^5, &  {\rm ord}(tr_5)=8, & ch_6=[t,tr_1]^2, & {\rm ord}(tr_6)=8,\\

ch_7=[t,tr_1]^3, & {\rm ord}(ch_7)=6, & ch_8=[t,tr_1]^4, & {\rm ord}(ch_6)=4.
\end{array}
$$

The set $\overline C_{7,\,5}$ consists of the following 10 ireducible semi-invariants:
$$
\begin{array}{llll}
pt_1=[t,ch_6]^5, &{\rm ord}(pt_1)=5, & pt_2=[t,ch_6]^6, & {\rm ord}(pt_2)=3,\\

pt_3=[t,ch_7]^2, & {\rm ord}(pt_3)=9, & pt_4=[t,ch_7]^3,  & {\rm ord}(pt_4)=7,\\

pt_5=[t,ch_7]^5,  & {\rm ord}(pt_5)=3, & pt_6=[t,ch_6]^3,  & {\rm ord}(pt_6)=9,\\

pt_7=[t,ch_4]^2,  & {\rm ord}(pt_7)=13, & pt_8=[t,ch_4]^5, & {\rm ord}(pt_6)=7, \\

pt_9=[t,dv^2_1]^7,  & {\rm ord}(pt_9)=5, & pt_{10}=[t,dv_1 dv_2]^7,  & {\rm ord}(pt_{10})=1. \\
\end{array}
$$

The set $\overline C_{7,\,6}$ consists of the following 10 ireducible semi-invariants:

$$
\begin{array}{llll}
sh_1=[t,pt_5]^5, & {\rm ord}(sh_1)=6, & sh_2=[t,pt_7]^6, &  {\rm ord}(sh_2)=8,\\

sh_3=[t,pt_4]^5, &  {\rm ord}(sh_3)=4, & sh_4=[t,pt_4]^6, & {\rm ord}(sh_4)=2,\\

sh_5=[t,pt_3]^2, &  {\rm ord}(sh_5)=12, & sh_6=[t,pt_3]^4,  & {\rm ord}(sh_6)=8,\\

sh_7=[t,pt_4]^4,  & {\rm ord}(sh_7)=6, & sh_8=[t,tr_1 dv_1]^7, & {\rm ord}(sh_6)=4, \\

sh_9=[t,tr_1 dv_2]^6,  & {\rm ord}(sh_9)=2, & sh_{10}=[t, tr_6 dv_1]^7, &  {\rm ord}(sh_{10})=2. \\
\end{array}
$$
       
 The set $\overline C_{7,\,7}$ consists of the following 12 ireducible semi-invariants: 

$$
\begin{array}{llll}
si_1=[t,sh_5]^4, & {\rm ord}(si_1)=11, & si_2=[t,sh_7]^4, & {\rm ord}(si_2)=5,\\

si_3=[t,tr_1^2]^7,  & {\rm ord}(si_3)=3, & si_4=[t,sh_1]^3, & {\rm ord}(si_4)=7,\\

si_5=[t,ch_7 dv_1]^7,  & {\rm ord}(si_5)=5, & si_6=[t,ch_7 dv_2]^7,  & {\rm ord}(si_6)=1,\\

si_7=[t,tr_6^2]^4, & {\rm ord}(si_7)=5, & si_8=[t,tr_6^2]^6,  & {\rm ord}(si_6)=1, \\

si_9=[t,tr_6\, tr_1]^6,  & {\rm ord}(si_9)=3, & si_{10}=[t, tr_6\, tr_1]^7, & {\rm ord}(si_{10})=1, \\

si_{11}=[t,tr_1^2]^6,  & {\rm ord}(si_{11})=5, & si_{12}=[t,sh_{10}), & {\rm ord}(si_{12})=7.
\end{array}
$$      

   The set $\overline C_{7,\,8}$ consists of the following 13 ireducible semi-invariants:     
$$
\begin{array}{llll}
vi_1=[t,si_7]^3, & {\rm ord}(vi_1)=6, & vi_2=[t,si_7]^4,  & {\rm ord}(vi_2)=4,\\

vi_3=[t,ch_8\,tr_6]^7,  & {\rm ord}(vi_3)=0, & vi_4=[t,ch_8\,tr_1]^6,  & {\rm ord}(vi_4)=4,\\

vi_5=[t,ch_8 tr_1]^7,  & {\rm ord}(vi_5)=2, & vi_6=[t,ch_7 tr_6]^7,  & {\rm ord}(vi_6)=2,\\

vi_7=[t,ch_7 \,tr_1]^7,  & {\rm ord}(vi_7)=4, & vi_8=[t,ch_8\,tr_6]^6, & {\rm ord}(vi_6)=2, \\

vi_9=[t,tr_6\, dv_2^2]^7, & {\rm ord}(vi_9)=0, & vi_{10}=[t,si_4]^2, & {\rm ord}(vi_{10})=10, \\

vi_{11}=[t,si_12]^4, & {\rm ord}(vi_{11})=6, & vi_{12}=[t,si_{11}]^3),  & {\rm ord}(vi_{12})=6,\\
& vi_{13}=[t,pt_9\,dv_2]^7, & {\rm ord}(vi_{13})=0.&
\end{array}
$$

The set $\overline C_{7,\,9}$ consists of the following 11 ireducible semi-invariants:
$$
\begin{array}{llll}
de_1=[t,sh_3\,dv_1]^7, & {\rm ord}(de_1)=3, & de_2=[t,ch_7\,ch_8]^7, & {\rm ord}(de_2)=3,\\

de_3=[t,pt_5\,tr_6]^5, & {\rm ord}(de_3)=3, & de_4=[t,pt_5 \,tr_1]^6,  & {\rm ord}(de_4)=3,\\

de_5=[t,pt_5\, tr_1]^7,  & {\rm ord}(de_5)=1, & de_6=[t,sh_9\,dv_1]^7,  & {\rm ord}(de_6)=1,\\

de_7=[t,sh_{10}\,dv_1]^7,  & {\rm ord}(de_7)=1, & de_8=[t,sh_{10} \,dv_2]^3, & {\rm ord}(de_6)=5, \\

de_9=[t,vi_5]^2,  & {\rm ord}(de_9)=5, & de_{10}=[t,vi_2]^4, & {\rm ord}(de_{10})=3, \\

& de_{11}=[t,vi_{11}]^2, & {\rm ord}(de_{11})=9.&
\end{array}
$$

The set $\overline C_{7,\,10}$ consists of the following 9 ireducible semi-invariants:
$$
\begin{array}{llll}
des_1=[t,sh_9\,tr_1]^6, & {\rm ord}(des_1)=2, & des_2=[t,sh_4\,tr_6]^4, & {\rm ord}(des_2)=4,\\

des_3=[t,sh_4\,tr_1]^6,  & {\rm ord}(des_3)=2, & des_4=[t,sh_1\,tr_1]^7, & {\rm ord}(des_4)=4,\\

des_5=[t,sh_3\,tr_6]^5, & {\rm ord}(des_5)=4, & des_6=[t,de_9]^2, & {\rm ord}(des_6)=8,\\

des_7=[t,tr_6^3]^7, & {\rm ord}(des_7)=2, & des_8=[t,sh_{10} \,tr_1]^6, & {\rm ord}(des_6)=2, \\

& des_9=[t,pt_1\,ch_7]^7,  & {\rm ord}(des_9)=4. &
\end{array}
$$

The set $\overline C_{7,\,11}$ consists of the following 9 ireducible semi-invariants:
$$
\begin{array}{llll}
odn_1=[t,vi_2\,dv_1]^7, &{\rm ord}(odn_1)=3, & odn_2=[t,vi_2,dv_2]^6,  & {\rm ord}(odn_2)=1,\\

odn_3=[t,vi_4\, dv_2]^6, & {\rm ord}(odn_3)=1, & odn_4=[t,vi_5\,dv_1]^7, & {\rm ord}(odn_4)=1,\\

odn_5=[t,vi_6\, dv_1]^7, & {\rm ord}(odn_5)=1, & odn_6=[t,vi_2\,dv_2]^5,  & {\rm ord}(odn_6)=3,\\

odn_7=[t,des_6]^4, & {\rm ord}(odn_7)=7, & odn_8=[t,des_6]^6,  &{\rm ord}(odn_6)=3, \\

& odn_9=[t,vi_1\,dv_2]^7, & {\rm ord}(odn_9)=1.&
\end{array}
$$

The set $\overline C_{7,\,12}$ consists of the following 13 ireducible semi-invariants:
$$
\begin{array}{llll}
dvan_1=[t,sh_1\,pt_2]^7, & {\rm ord}(dvan_1)=2, & dvan_2=[t,sh_1\,pt_5]^7, & {\rm ord}(dvan_2)=2,\\

dvan_3=[sh_9,sh_{10}]^2, & {\rm ord}(dvan_3)=0, & dvan_4=[t,odn_7]^6, & {\rm ord}(dvan_4)=2,\\

dvan_5=[t, de_8\,dv_2]^6, & {\rm ord}(dvan_5)=2, & dvan_6=[sh_{10}\, ,sh_{10}]^2, & {\rm ord}(dvan_6)=0,\\

dvan_7=[t, de_9\,dv_2]^6, & {\rm ord}(dvan_7)=2, & dvan_8=[t, de_{10}\,dv_1]^7,  & {\rm ord}(dvan_6)=2, \\

dvan_9=[t,odn_7]^4, & {\rm ord}(dvan_9)=6, & dvan_{10}=[sh_1\, , sh_1]^2, & {\rm ord}(dvan_{10})=0, \\

dvan_{11}=[sh_4\, , sh_4]^2, & {\rm ord}(dvan_{11})=0, & dvan_{12}=[sh_4, sh_9]^2), & {\rm ord}(dvan_{12})=0,\\

&  dvan_{13}=[sh_4\, , sh_2]^2, & {\rm ord}(dvan_{13})=0. &
\end{array}
$$   
  
The set $\overline C_{7,\,13}$ consists of the following 9 ireducible semi-invariants:
$$
\begin{array}{llll}
tryn_1=[t,dvan_9]^6, & {\rm ord}(tryn_1)=1, & tryn_2=[t,vi_1\,ch_7]^7, &  {\rm ord}(tryn_2)=5,\\

tryn_3=[t,vi_2\,ch_8]^7, & {\rm ord}(tryn_3)=1, & tryn_4=[t,vi_2\,ch_2]^7, & {\rm ord}(tryn_4)=1,\\

tryn_5=[t,vi_1\,ch_8]^7, & {\rm ord}(tryn_5)=3, & tryn_6=[t,vi_5\,ch_2]^6, & {\rm ord}(tryn_6)=1,\\

tryn_7=[t,vi_8\,ch_8]^6, & {\rm ord}(tryn_7)=1, & tryn_8=[t,vi_8\,ch_7]^7,  & {\rm ord}(tryn_6)=1, \\

& tryn_9=[t,vi_4\,ch_8]^7, & {\rm ord}(tryn_9)=1. &

\end{array}
$$  

A number of elements of  $\overline C_{7,\,i},$ $i=1,\ldots,13$  and their orders so far coincide completelly with Gall's  results, see   \cite{Gall-1}.

%----------------------------------------------------------------------------------

  \section{Computation of the invariants.}

%----------------------------------------------------------------------------------------------------

Put  $I_{i}:=\bar I_7 \cap \overline C_{7,\,i},$ $I_{+}:=I_7^{+}.$ Let  $(I_{+}^2)_i$ be a subset  of  $(I_7^{+})^2$ whose elements  has degree  $i.$ 
Denote by  $\delta_i$ a  number of linearly independent irreducible invariants of degree $i.$  It is evident that    ${\delta_i =\dim I_i-\dim(I^2_{+})_{i}.}$ A dimension of the vector space  $I_i$ is calculating by Cayley-Sylvester formula, see, for example  \cite{Hilb}. A dimension of the vector space  $(I^2_{+})_i$ is calculating by the formula $\dim(I^2_{+})_i=\sigma_i -\dim S_i.$ Here  $\sigma_i$ is coefficient of  $x^i$ in the series expansion  $\displaystyle \left(\prod_{k<i} (1-x^k)^{\delta_k}\right)^{-1},$ and  $S_i$ is a vector space of  $(I^2_{+})_i$ generated by syzygies. 

The invariants of  4, 8 and  12 degries were found above.  We have
$$
\begin{array}{l}
I_4=\langle p_4  \rangle, \delta_4=1, p_4:=ch_1=[t,tr_5]^7,\\
I_8=\langle p_{8,1},p_{8,2},p_{8,3}\rangle , \delta_8=3
 \\p_{8,1}:=vi_3=[t,ch_8\,tr_6]^7,p_{8,2}:=vi_9=[t,tr_6\, dv_2^2]^7,p_{8,3}:=vi_{13}=[t,pt_9\,dv_2]^7,\\

I_{12}=\langle p_{12,1},p_{12,2},p_{12,3},p_{12,4},p_{12,5},p_{12,6}\rangle, \delta_{12}=6 \\

p_{12,1}:=dvan_3=[sh_9,sh_{10}]^2, p_{12,2}:=dvan_6=[sh_{10}\,sh_{10}]^2,p_{12,3}:= dvan_{10}=[sh_1\, sh_1]^2, \\

p_{12,4}:=dvan_{11}=[sh_4\,sh_4]^2, p_{12,5}:= dvan_{12}=[sh_4,sh_9]^2),p_{12,6}:=dvan_{13}=[sh_4\,sh_2]^2.
\end{array}
$$

For $I_{14}$ we have  $\dim I_{14}=4,$ $\sigma_{14}=0.$  Therefore  $(I_{+}^2)_{14}=0,$ then  $\delta_{14}=4.$ It is enought to find  $4$ linearly independent invariants of degree $14.$ Below is typical instance of how these  invariants are calculated.  Invariants of $I_{14}$ we seek to find are semitransvectants of the form $[u,v]^i,$ where  $u, v$ are semi-invariants of $\overline C_{7,7},$ $i<4,$ ${\rm ord}(u)+{\rm ord}(v) - 2\,i=0.$ There are  $9$ such semitransvectants:
$$
[si_8,si_8], [si_8,si_{10}], [si_9,si_9]^3, [si_{10},si_{10}], [si_6,si_6], [si_6,si_8], [si_6,si_{10}], [si_3, si_3]^3, [si_3,si_9]^3.
$$
By using Maple calculation we choose four  linear independent invariants -- $$p_{14,1}:=[si_8,si_{10}], p_{14,2}:= [si_6,si_{10}], p_{14,3}:=[si_6,si_8], p_{14,4}:=[si_3,si_9]^3.$$
The invariants  $p_{14,1}, p_{14,2}, p_{14,3},p_{14,4}$ are fractions with the denominator  $t^{35}.$ The numerators of the fractions are polynomials of  $\mathbb{Z}[z_2,\ldots,z_7]$ which consists of    937,  869,  978,   925  terms respectivelly.

For $I_{16}$ we have  $\dim I_{16}=18,$ $\sigma_{16}=16.$ The vector space  $(I_{+}^2)_{16}$ is generated by 16 elements and all of them are linearly independent.  Thus $\delta_{16}=2.$ In order to define two invariants of   $I_{16}$ let us consider a set of semitransvectants of the form  $[u,v]^i,$ where  $u, v \in \overline C_{7,8},$ $i<5.$ There are $12$  such semitransvectants
$$
\begin{array}{l}

[vi_7,vi_7]^4, [vi_5,vi_5]^2, [vi_5,vi_6]^2, [vi_5,vi_8]^2, [vi_6,vi_6]^2, [vi_6,vi_8]^2,\\

 [vi_2,vi_7]^4, [vi_4,vi_7]^4, [vi_2,vi_2]^4, [vi_2,vi_4]^4, [vi_8,vi_8]^2, [vi_4,vi_4]^4.

\end{array}
$$

In order to separate 2 linearly independent irreducible semitransvectants let us consider the equality: 

$$
\alpha_1 p_4^4+\alpha_2 p_8^2+\ldots +\alpha_{17} [vi_7,vi_7]^4+\ldots +\alpha_{28} [vi_4,vi_4]^4=0.
$$
Substituting the invariants values into the equality we obtain an overdefined system of equations. After solving the system we get that the following 18 elements --  16 basis elements of the vectors space $(I_{+}^2)_{16}$ and the two invariants $p_{16,1}:=[vi_2,vi_4]^4,$ and ${p_{16,2}:=[vi_4,vi_7]^4}$  span the vector space $I_{16}.$ 

The invariants  $p_{16,1}, p_{16,2}$ are fractions with the denominator  $t^{40}.$ The numerators of the fractions are  polynomials of  $\mathbb{Z}[z_2,\ldots,z_7]$ which consist of   1744 and   1698  terms.

For $I_{18}$ we have  $\dim I_{18}=13,$ $\sigma_{18}=4.$ Since  $(I_{+}^2)_{18}=p_4 I_{14},$ then  $S_{18}=0.$ Thus  $\delta_{18}=9.$  
Invariants of $I_{18}$ we seek to find are semitransvectants of the form $[u,v]^i,$ where  $u, v$ are semi-invariants of $\overline C_{7,9},$ $i<6,$ ${\rm ord}(u)+{\rm ord}(v) - 2\,i=0$ 

In the same way as above we obtain the nine irreducible invariants : 
$$
\begin{array}{lll}
p_{18,1}:=[de_4,de_3]^ 3, & p_{18,2}:=[de_4,de_{10}]^ 3, & p_{18,3}:=[de_5,de_6], \\
p_{18,4}:=[de_1,de_{10}]^ 3, & p_{18,5}:=[de_2,de_3]^ 3, & p_{18,6}:=[de_2,de_{10}]^ 3,\\
 p_{18,7}:=[de_3,de_{10}]^ 3, & p_{18,8}:=[de_6,de_{7}],  & p_{18,9}:=[de_8,de_{9}]^5.
\end{array}
$$
The invariants  $p_{18,1}, p_{18,2}, \ldots, p_{18,9}$ are fractions with the denominator  $t^{45}.$ The numerators of the fractions are  polynomials of  $\mathbb{Z}[z_2,\ldots,z_7]$ which consist of   2674, 2758, 2645, 2800, 2718, 2772, 2769, 2661, 2739  terms respectivelly.

For $I_{20}$  we have  $\dim I_{20}=35,$ $\sigma_{20}=36.$ The vector space $S_{20}$ is spaned  by two syzygies: 
$$
\begin{array}{l}
 - 142725300\,{p_{12, \,4}}\,{p_{4}}^{2} - 15449224200\,{p_{8, \,
1}}\,{p_{8, \,2}}\,{p_{4}} + 1320855600\,{p_{8, \,3}}\,{p_{8, \,1
}}\,{p_{4}}+\\

\mbox{} - 327281580\,{p_{8, \,3}}\,{p_{8, \,2}}\,{p_{4}} - 
75375000\,{p_{12, \,6}}\,{p_{8, \,1}} + 13989210\,{p_{8, \,3}}^{2
}\,{p_{4}}  - 691200\,{p_{16, \,2}}\,{p_{4}}+ \\

 {+} 1530000\,{p_{12, \,4}}\,{p_{8, \,3}}{-}19980000\,{p_{12
, \,5}}\,{p_{8, \,2}} {+} 31075755000\,{p_{8, \,1}}^{2}\,{p_{4}}{+}1890487920\,{p_{8, \,2}}^{2}\,{p_{4}}{+}

 \\

{-}16875000\,{p_{12, \,2}}\,{p_{8, \,1}}{-}2025000\,{p_{12
, \,1}}\,{p_{8, \,3}} + 34290000\,{p_{12, \,1}}\,{p_{8, \,2}}{-}151200000\,{p_{12, \,1}}\,{p_{8, \,1}} +
\\

\mbox{} + 360450000\,{p_{12, \,5}}\,{p_{8, \,1}} - 10800000\,{p_{
12, \,4}}\,{p_{8, \,2}} + 197100000\,{p_{12, \,4}}\,{p_{8, \,1}} - 675000\,{p_{16
, \,1}}\,{p_{4}} +
 \\

\mbox{} - 337500\,{p_{12, \,2}}\,{p_{8, \,3}} - 1507500\,{p_{12, 
\,6}}\,{p_{8, \,3}} + 2970000\,{p_{12, \,5}}\,{p_{8, \,3}} - 188556795\,{p_{12, \,5}}\,{p_{4}}^{2}+  \\

\mbox{} + 4050000\,{p_{12, \,2}}\,{p_{8, \,2}} + 19440000\,{p_{12
, \,6}}\,{p_{8, \,2}} - 135229030\,{p_{8, \,1}}\,{p_{4}}^{3}+ 24000
\,{p_{12, \,2}}\,{p_{4}}^{2} + \\
\mbox{}+ 35439332\,{p_{
8, \,2}}\,{p_{4}}^{3} - 2797842\,{p_{8, \,3}}\,{p_{4}}^{3} +110336985\,{p_{12, \,1
}}\,{p_{4}}^{2}   =0,

\end{array}
$$
and  
$$
\begin{array}{l}
144093300\,{p_{12, \,4}}\,{p_{4}}^{2} + 15499144200\,{p_{8, \,1}}
\,{p_{8, \,2}}\,{p_{4}} - 1328535600\,{p_{8, \,3}}\,{p_{8, \,1}}
\,{p_{4}}  +\\

{+}329009580\,{p_{8, \,3}}\,{p_{8, \,2}}\,{p_{4}}{-} 
68625000\,{p_{12, \,6}}\,{p_{8, \,1}}{-}14085210\,{p_{8, \,3}}^{2
}\,{p_{4}}{-}110768985\,{p_{12, \,1}}\,{p_{4}}^{2}{-}\\

{-} 4410000\,{p_{12, \,4}}\,{p_{8, \,3}}{-}49140000\,{p_{12
, \,5}}\,{p_{8, \,2}}{-}31229355000\,{p_{8, \,1}}^{2}\,{p_{4}}{-}1903639920\,{p_{8, \,2}}^{2}\,{p_{4}}{+}
 \\

\mbox{} + 16875000\,{p_{12, \,2}}\,{p_{8, \,1}} + 2025000\,{p_{12
, \,1}}\,{p_{8, \,3}} + 270000\,{p_{12, \,1}}\,{p_{8, \,2}} + 675000\,{p_{16
, \,1}}\,{p_{4}}+\\

\mbox{} + 676350000\,{p_{12, \,5}}\,{p_{8, \,1}} + 10800000\,{p_{
12, \,4}}\,{p_{8, \,2}} - 110700000\,{p_{12, \,4}}\,{p_{8, \,1}}+
 \\

\mbox{} + 337500\,{p_{12, \,2}}\,{p_{8, \,3}} - 1372500\,{p_{12, 
\,6}}\,{p_{8, \,3}} + 14310000\,{p_{12, \,5}}\,{p_{8, \,3}}  + 691200
\,{p_{16, \,2}}\,{p_{4}}-\\

\mbox{} - 4050000\,{p_{12, \,2}}\,{p_{8, \,2}} + 15120000\,{p_{12
, \,6}}\,{p_{8, \,2}} + 136157030\,{p_{8, \,1}}\,{p_{4}}^{3}+ 336000\,{p_{12
, \,6}}\,{p_{4}}^{2}- \\
\mbox{} - 21600000\,{p_{12, \,1}}\,{p_{8, \,1}} - 35682532\,{p_{8
, \,2}}\,{p_{4}}^{3} + 2817042\,{p_{8, \,3}}\,{p_{4}}^{3}+ 186396795\,{p_{12, \,5}}\,{p_{4}}^{2}=0. \\

\end{array}
$$
Hence  $\dim S_{20}=2$ and we get    $\delta_{20}=1.$ It confirm the results of the paper \cite{DL}.

 The unique invariant $I_{20}$ we are considering is a semitransvectants of the form  $[u,v]^i,$ where $u, v \in \overline C_{7,10},$ $i<3,$ ${\rm ord}(u)+{\rm ord}(v) - 2\,i=0.$  As result of a calculation we get that the element  $p_{20}:=[des_7,des_7]^2$ is irreducle invariant of degree 20. 
The invariant  $p_{20}$ is  a fraction with the denominator  $t^{50}.$ The numerator of the fraction is  a polynomial of  $\mathbb{Z}[z_2,\ldots,z_7]$ which consists of   4392  terms.

For $I_{22}$ we have  $\dim I_{22}=26,$ $\sigma_{22}=25.$ By direct calculation we  found that those 25  elements of  $(I_{+}^2)_{22}$ satisfy the unique syzygy:
$$
\begin{array}{l}
{p_{14, \,3}}\,{p_{8, \,3}} - 4\,{p_{14, \,2}}\,{p_{8, \,3}} + 40
\,{p_{14, \,3}}\,{p_{8, \,1}} + 750\,{p_{14, \,1}}\,{p_{8, \,3}}
 + 150\,{p_{14, \,4}}\,{p_{8, \,3}} - 160\,{p_{14, \,2}}\,{p_{8, 
\,1}} -\\
\mbox{} - 1275\,{p_{14, \,4}}\,{p_{8, \,2}} - 21\,{p_{14, \,3}}\,
{p_{8, \,2}} - 1875\,{p_{18, \,8}}\,{p_{4}} + 84\,{p_{14, \,2}}\,
{p_{8, \,2}} + 48750\,{p_{14, \,1}}\,{p_{8, \,1}}+ \\
\mbox{} + 9750\,{p_{14, \,4}}\,{p_{8, \,1}} + 1125\,{p_{14, \,1}}
\,{p_{8, \,2}}=0.
\end{array}
$$
Hence $\dim S_{22}=1$ and    $\delta_{22}=26-25+1=2.$ It coincides with the results of  \cite{DL}.

 Invariants of $I_{22}$ we are looking as semitransvectants of the form  $[u,v]^i,$ where   $u, v \in \overline C_{7,11},$ $i<4,$ ${\rm ord}(u)+{\rm ord}(v) - 2\,i=0.$  As result of the  calculation we get that the elements  $$p_{22,1}:=[odn_6, odn_1]^3, p_{22,2}:=[odn_8, odn_1]^3$$ are linearly independent  irreducible invariants of degree 22. 

The invariants  $p_{22,1},$ $p_{22,2}$ are fractions with the denominator  $t^{55}.$ The numerators of the fractions are  polynomials of  $\mathbb{Z}[z_2,\ldots,z_7]$ which consist of   6 569 and  6 556  terms respectivelly.

For $I_{24}$ we have  $\dim I_{24}=62,$ $\sigma_{24}=74.$ By direct calculation we get the  vector space $S_{24}$ is spaned by 12 syzygies. Hence  $\delta_{24}=62-74+12=0. $

For $I_{26}$ we have  $\dim I_{26}=52,$ $\sigma_{26}=78.$ By direct Maple calculation, we obtain  that the  vector space $S_{26}$ is spaned by 27 syzygies. Hence, $\delta_{26}=52-78+27=1, $ that coincides with the results of  \cite{DL}.
Invariants of $I_{26}$ we are looking as semitransvectants of the form  $[u,v]^i,$ where   $u, v \in \overline C_{7,13},$ $i=1,$ ${\rm ord}(u)+{\rm ord}(v) - 2\,i=0. $
As result of the  calculation we get that the element $$p_{26}=[tryn_4,tryn_3]$$ is the irreducible, unique invariant of degree 26. 
The invariants  $p_{26}$ is a fraction with the denominator  $t^{65}.$ The numerator are  polynomials of  $\mathbb{Z}[z_2,\ldots,z_7]$ which consist of  13 651  terms.
  
For $I_{28}$ we have  $\dim I_{28}=97,$ $\sigma_{28}=135.$ By direct calculation we get the  vector space $S_{28}$ is spaned by 38 syzygies. Hence  $\delta_{28}=97-135+38=0. $

For $I_{30}$ we have  $\dim I_{30}=92,$ $\sigma_{30}=171.$ 
By direct Maple calculation we obtain  that the  vector space $S_{30}$ is spaned by 80 syzygies. Hence, $\delta_{30}=92-171+80=1.$ 

The unique irreducible invariant of $I_{30}$ we take from the paper   \cite{Gall-1} where the notation  ${p_{30}=(h,\alpha)}$  is used. After all the calculations,  we obtain that  $p_{30}$ is fraction with the denominator  $t^{75}.$ The numerator of the fraction are polynomials of  $\mathbb{Z}[z_2,\ldots,z_7]$ which consist of  25 868  terms.

 Summarizing the above results we get
\begin{te}
The  system of the 30 invariants
$$
\begin{array}{l}
 p_4=[t,tr_5]^7,\\
\\
 p_{8\!,1}=[t,ch_8\,tr_6]^7,p_{8\!,2}=[t,tr_6\, dv_2^2]^7,p_{8,3}=[t,pt_9\,dv_2]^7,\\
\\
p_{12,1}=[sh_9,sh_{10}]^2, p_{12,2}=[sh_{10}\,sh_{10}]^2,p_{12,3}=[sh_1\, sh_1]^2, \\
p_{12,4}=[sh_4\,sh_4]^2, p_{12,5}=[sh_4,sh_9]^2),p_{12,6}=[sh_4\,sh_2]^2,\\
\\
p_{14,1}=[si_8,si_{10}], p_{14,2}= [si_6,si_{10}], p_{14,3}=[si_6,si_8], p_{14,4}=[si_3,si_9]^3,\\
\\
p_{16,1}=[vi_2,vi_4]^4, p_{16,2}=[vi_4,vi_7]^4,\\
\\
p_{18,1}=[de_4,de_3]^ 3, p_{18,2}=[de_4,de_{10}]^ 3, p_{18,3}=[de_5,de_6], p_{18,4}=[de_1,de_{10}]^ 3, p_{18,5}=[de_2,de_3]^ 3, \\

p_{18,6}=[de_2,de_{10}]^ 3, p_{18,7}=[de_3,de_{10}]^ 3, p_{18,8}=[de_6,de_{7}],  p_{18,9}=[de_8,de_{9}]^5,\\
\\
p_{20}=[des_7,des_7]^2,\\
\\
p_{22,1}=[odn_6, odn_1]^3, p_{22,2}=[odn_8, odn_1]^3,\\
\\
p_{26}=[tryn_4,tryn_3],\\
\\
p_{30}=(h,\alpha) \mbox{{ (in Gall's notation.)}}
\end{array}
$$

 is a complete system of the invariants for the  binary form of  degree 7.
\end{te}
  The explicit  form of these invariants  is highly unwieldy: a text  file listing these invariants is as large as 4MB.

\end{document}